\documentclass[9pt,english]{article}
\usepackage{geometry}
\usepackage{amsmath,amsfonts, amsthm,amssymb,oldgerm,amssymb}
\usepackage{array}
\usepackage[english]{babel}
\usepackage{float}
\usepackage{url}
\usepackage{cite} 
\usepackage{hyperref}
\usepackage{enumitem}
\usepackage{xcolor} 
\usepackage{amssymb,mathrsfs,graphicx,amsmath,graphicx,amsthm}
 \numberwithin{equation}{subsection}
\usepackage{mathtools}
\usepackage{array,multicol,comment,lipsum,caption}%bibtex}
\usepackage{color}

\theoremstyle{plain}
\newtheorem{theorem}{Theorem}[section]
\newtheorem{proposition}[theorem]{Proposition}

\theoremstyle{definition}
\newtheorem{definition}{Definition}[section]

\theoremstyle{plain}
\newtheorem{lemma}[theorem]{Lemma}

\theoremstyle{remark}
\newtheorem{remark}{Remark}[section]

\usepackage{epsfig}
\usepackage{relsize}
\usepackage{hyperref}
\usepackage[utf8]{inputenc}

\newcommand{\Int}{\displaystyle\int}
\newcommand{\R}{{\mathbb{R}}}

\DeclarePairedDelimiter\norm{\lVert}{\rVert}
\DeclarePairedDelimiterX{\inp}[2]{\langle}{\rangle}{#1, #2}

\makeatletter
\newcommand{\myitem}[1]{%
\item[#1]\protected@edef\@currentlabel{#1}%
}
\makeatother

\newcommand{\mycomment}[1]{}
\begin{document}
\title{\bf Continuity and Stability of the FitzHugh-Nagumo System with Respect to the Diffusion}

\author{Alexander Rodr\'iguez\ \footnotemark[2]}

%\footnotetext[1]{This work has been partially funded by NNN}

\footnotetext[2]{Facultad de Matem\'aticas, Pontificia Universidad Católica de Chile, Avda.
Vicu\~na Mackenna 4860, Macul, Santiago, Chile. E.mail: jarodriguez19@uc.cl}

%\author{Your Name}
\date{\today}
\maketitle
\begin{abstract}
Spatial heterogeneity in diffusion is a fundamental feature of excitable media, yet its effect on FitzHugh–Nagumo dynamics remains mathematically challenging. We establish stability of solutions with respect to perturbations of the diffusion coefficient, providing a rigorous account of how changes in the medium affect the evolution of the system. Specifically, we prove a continuous dependence result ensuring that small perturbations in diffusion produce small variations in the corresponding solutions. The proof adapts the method of contracting rectangles to variable diffusion, yielding a framework for controlling and comparing solutions associated with different diffusion profiles. This result provides a rigorous link between heterogeneous and homogeneous FitzHugh–Nagumo systems and justifies constant-diffusion approximations when the diffusion coefficient is close to a constant. Our analysis advances the understanding of FitzHugh–Nagumo dynamics in heterogeneous media by establishing robustness with respect to spatial variations in diffusion.
\end{abstract}
\section{Introduction}
The propagation of electrical signals in biological tissues, such as neurons and cardiac tissue, is governed by reaction-diffusion systems that describe the interaction between excitation and recovery variables. The FitzHugh-Nagumo (FHN) model is a simplified version of the Hodgkin-Huxley equations that captures the essential dynamics of excitable media \cite{fhn,4066548nagumo,hodgkin1952quantitative,rattay1986analysis,MR1952568,MR2263523,MR381744}. From a mathematical perspective, the FHN system has served as a prototype for studying nonlinear parabolic systems, particularly in the context of traveling wave solutions, pulse propagation, and pattern formation \cite{rinzel1973traveling,MR889639KEENER,MR393759HASTING}. In its classical form, the FHN model assumes constant diffusion coefficients; however, in realistic applications, the diffusion is often spatially varying due to biological heterogeneities, which represents the tissue conductivity governing the spread of electrical activity.\\

In neuroscience, the conductivity of neural tissue depends on structural features such as axon diameter, myelination, and extracellular resistance \cite{rall1969time,MR889639KEENER}. These properties critically determine the speed and efficiency of signal transmission. The diffusion coefficient in the FHN model reflects this effective conductivity, influencing how excitation propagates through space. For example, larger axons and saltatory conduction enhance signal spread, while demyelination in diseases like multiple sclerosis reduces conductivity and impairs function \cite{rall1969time}. In cardiac tissue, similar effects occur: gap junction density, fibrosis, and ischemia all alter tissue conductivity. Fibrosis creates non-conductive regions that block wave propagation, while ischemia leads to heterogeneous conduction by disrupting ion channel function \cite{clayton2011models,kleber2004basic}. Thus, spatial variation in the diffusion coefficient provides a mathematical representation of structural inhomogeneities in biological tissues.\\

Mathematically, the FHN system is composed by two coupled parabolic equations and in this paper we consider it with a variable diffusion. Thus, the systems read as
\begin{equation}\label{delta-fhn}
\begin{aligned}
v_t & =\left(\delta(x) v_x\right)_x+f(v)-u\\
u_t & =\rho u_{xx} + \varepsilon(v-\gamma u)
\end{aligned}
\end{equation}
where $v,u$ are real-valued functions on the variables $x \in \mathbb{R}$ and $t \geq 0$, $\delta(x) > 0$ for all $x\in \R$, $\rho \geq 0$, the parameters $\varepsilon$ and $\gamma$ are chosen to be positive and $f$, the reaction term, is given by the cubic polynomial
\begin{equation}
f(v)=v(v-a)(1-v), \quad \text{with } a \in (0,1/2).
\end{equation}

In the case of constant diffusion, numerous results have established the well-posedness of the FitzHugh–Nagumo system, ensuring existence, uniqueness, and regularity of solutions under standard assumptions \cite{RAUCH197812,MR425397,SCHONBEK1978119,MR4583896,MR3927583}. A major line of research has focused on the existence and stability of traveling wave solutions, using tools such as geometric singular perturbation theory, the Maslov index, and spectral analysis \cite{MR3551276carterrijk,MR3395129cartersandstede,MR3773759cornwell,MR3281845arioli,MR3385151chen2,MR4419132sanchez}. The effects of diffusion and spatial heterogeneity have also been widely studied, as they play a crucial role in shaping solution behavior, wave propagation, and pattern formation \cite{MR4860996ammatouch,MR716853collins,MR2494028faina,MR4484075gambino}. On the computational side, traveling waves have been numerically approximated using, for instance, equilibria of associated nonlocal equations \cite{MR4887572rubio}. More recently, extended formulations have incorporated fractional-order derivatives and stochastic effects, allowing for the analysis and optimal control of FitzHugh–Nagumo-type dynamics in more general and realistic settings \cite{MR4274265orkun,MR4308218cordoni,cebrian2024sixgalens}.\\

When $\delta(x)$ varies smoothly around a constant $\sigma$, one expects that the solutions of the FHN model with variable diffusion remain close to the solutions of the constant-diffusion equation:
\begin{equation}\label{fhn-knwon}
\begin{aligned}
& \tilde{v}_t=\sigma\tilde{v}_{x x}+f(\tilde{v})-\tilde{u} \\
& \tilde{u}_t=\rho \tilde{u}_{xx}+\varepsilon(\tilde{v}-\gamma \tilde{u}).
\end{aligned}
\end{equation}

This works aims to establish the continuity of the solution map with respect to the diffusion coefficient in the FitzHugh-Nagumo model, proving that small variations in $\delta(x)$ lead to small variations in wave propagation. This stability result ensures that minor perturbations in tissue conductivity do not cause drastic changes in the solution, providing a rigorous justification for the use of effective homogeneous models to approximate excitable media when diffusion variations are small. The main idea of the work is inspired by the technique of contracting rectangles \cite{RAUCH197812}, which has been successfully applied to the analysis of FitzHugh-Nagumo systems in various settings. This approach has been successfully applied to study the FHN system in various contexts, including high-frequency forcing, boundary control, and perturbative regimes \cite{cerpa2023approximation,cerpaode,SCHONBEK1978119,RAUCH197812,MR445091,MR1036209,MR3262177, jhoncerpa}. By adapting this framework, we extend these ideas to the case of variable diffusion coefficients, further demonstrating the robustness of the FitzHugh-Nagumo system under spatial heterogeneity.

\section{Main result}
We start this section by introducing  the concept of a contracting rectangle and all the functional spaces appearing in this paper.
\begin{definition}\label{definition-contracting-rectangle}
\textbf{(Contracting Rectangle)}. Let $R\subset \mathbb{R}^2$ be a rectangle such that the $(0,0)$ is in the interior. Given a vector-valued function $H: \mathbb{R}^2 \times \R \times [0,T] \rightarrow \mathbb{R}^2$ we say that the rectangle $R$ is contracting under the vector fields $H$ if for each $\vec{v} \in \partial R$ and every outward normal unit vector $n(\vec{v})$ we have
$$
\sup_{x\in\R} n(\vec{v}) \cdot H(\vec{v}, x, t)<0, \quad t\in [0,T].
$$    
\end{definition}
\begin{definition}
We set the following spaces involving functions $U: \R \to \R$, for $k\geq 0$ an integer number:
$$
\begin{aligned}
&B C^k(\R)=\{U:(d / d x)^j U  \, \text{is a bounded uniformly continuous function on } \, \mathbb{R} \, \text{for} \, \, 0 \leqslant j \leqslant k\}, \\
& C_0^k(\R)=\left\{U \in B C^k(\R):\lim (d / d x)^j U=0 \text { as }|x| \rightarrow \infty \text { for } 0 \leqslant j \leqslant k\right\},\\
&W^{k,p}(\R)=\left\{U \in L^p\left(\mathbb{R}\right):(d / d x)^j U \in L^p\left(\mathbb{R}\right) \text { for } 0 \leqslant j \leqslant k\right\}.
\end{aligned}
$$
%The Sobolev space $W^{k, p}$ consists of functions whose classical derivatives up to order $k$ exist and belong to $L^p$.
\end{definition}
We also define the following
\begin{definition}
For any $f \in C\left([0, T] ; B C^0(\R)\right)$ we consider the norm
$$
\|f\|_Y=\sup _{0<t<T} \sup _{x \in \mathbb{R}}|f(t,x)| .
$$
Additionally, for a fixed $t$, we denote
$$
\|f(t,\cdot)\|_{L_x^{\infty}}=\sup _{x \in \mathbb{R^+}}|f(t,x)| .
$$
\end{definition}
Concerning the functional framework where our results will hold.
\begin{definition}
Let $\mathcal{S}(t)$ denote the linear evolution operators
associated with the diffusion system. A Banach space
$\mathcal{B}$ of functions $w:\mathbb{R}\to\mathbb{R}$
is admissible for the diffusion system and the nonlinearity
$F$ if the following conditions hold. We equip
$\mathcal{B}^n$ with the norm
\[
\|U\|_{\mathcal{B}^n}
=\max_{1\leq i\leq n}\|u_i\|_{\mathcal{B}}.
\]
\begin{enumerate}
\myitem{A.1}\label{A.1}
The space $\mathcal{B}$ is continuously embedded in
$BC^0(\mathbb{R})$, the space of bounded uniformly
continuous functions, and
\[
\|w\|_{L^\infty(\mathbb{R})}
\leq\|w\|_{\mathcal{B}},
\qquad w\in\mathcal{B}.
\]

\myitem{A.2}\label{A.2}
The operators $\mathcal{S}(t)$ preserve $\mathcal{B}^n$
and form a strongly continuous semigroup of bounded linear
operators on $\mathcal{B}^n$. That is,
\[
\mathcal{S}(0)=I,
\qquad
\mathcal{S}(t+s)=\mathcal{S}(t)\mathcal{S}(s),
\qquad t,s\geq0,
\]
and
\[
\lim_{t\downarrow0}
\|\mathcal{S}(t)U-U\|_{\mathcal{B}^n}=0,
\qquad U\in\mathcal{B}^n.
\]
Consequently, there exists $K\geq1$ such that
\[
\sup_{0\leq t\leq1}
\|\mathcal{S}(t)\|_{\mathcal{L}(\mathcal{B}^n)}
\leq K.
\]

\myitem{A.3}\label{A.4}
The nonlinearity $F:\mathbb{R}^n\to\mathbb{R}^n$
satisfies $F(0)=0$, and its composition operator maps
$\mathcal{B}^n$ into itself. For every $M>0$, there exist
constants $k_1,k_2>0$ such that, for all
$V,W\in\mathcal{B}^n$,
\[
\begin{aligned}
\|V\|_{\mathcal{B}^n},\|W\|_{\mathcal{B}^n}\leq M
&\quad\Longrightarrow\quad
\|F(V)-F(W)\|_{\mathcal{B}^n}
\leq k_1\|V-W\|_{\mathcal{B}^n},\\
\|V\|_{L^\infty(\mathbb{R};\mathbb{R}^n)}\leq M
&\quad\Longrightarrow\quad
\|F(V)\|_{\mathcal{B}^n}
\leq k_2\|V\|_{\mathcal{B}^n}.
\end{aligned}
\]
Here the $L^\infty$ norm of a vector-valued function
is defined using the maximum norm on $\mathbb{R}^n$.
\end{enumerate}
\end{definition}

We now present our two main results. The first one establishes the existence and uniform boundedness of solutions to the FitzHugh–Nagumo system with variable diffusion, provided the diffusion coefficient remains close to a constant. The second one shows that, under suitable assumptions, the solution of the variable-diffusion system can be approximated by the solution of the corresponding constant-diffusion system.

\begin{theorem}[Global existence and exponential decay under variable diffusion]
\label{prop-boundedness}
Let $\mathcal B$ continuously embedded in
$C_0(\mathbb R)$ and satisfying assumptions \ref{A.1}--\ref{A.4}, let $a\in(0,1/2)$,
and let $\gamma>0$ such that $a>\frac{1}{\gamma}$. Let $\delta\in C^1(\mathbb{R})$ be a variable
diffusion coefficient such that $\delta$ and its derivative $\delta'$
are globally H\"older continuous with some exponent $\theta\in(0,1)$.
Assume that there exist a constant $\delta_*\geq1$ such that
\[
\delta_*^{-1}\leq\delta(x)\leq\delta_*,
\qquad x\in\mathbb{R}.
\] 
Then, for every initial data
$(v_0,u_0)\in\mathcal{B}\times\mathcal{B}$, there exists a contracting
rectangle $R\subset\mathbb{R}^2$ with
\[
(v_0(x),u_0(x))\in R
\qquad \text{for all } x\in\mathbb{R},
\]
such that the variable-diffusion system \eqref{delta-fhn} admits a
unique global solution
\[
(v,u)\in C\bigl([0,\infty),\mathcal{B}\times\mathcal{B}\bigr)
\]
satisfying
\begin{equation}
(v(x,t),u(x,t))\in R,
\qquad \forall x\in\mathbb{R},\ t\geq0.
\end{equation}
Moreover, there exist constants $\eta>0$, $C\geq 1$, and $c>0$,
depending only on $a$, $\gamma$, and $\varepsilon$, such that, if
\[
\max\bigl\{\|v_0\|_{L^\infty(\mathbb R)},
           \|u_0\|_{L^\infty(\mathbb R)}\bigr\}\leq \eta,
\]
then the solution satisfies
\[
\max\bigl\{\|v(t)\|_{L^\infty(\mathbb R)},
           \|u(t)\|_{L^\infty(\mathbb R)}\bigr\}
\leq
C e^{-ct}
\max\bigl\{\|v_0\|_{L^\infty(\mathbb R)},
           \|u_0\|_{L^\infty(\mathbb R)}\bigr\},
\qquad t\geq 0.
\]
\end{theorem}

\begin{theorem}[Uniform approximation by the constant-diffusion FHN system]
\label{thm-approximation}
Let $\mathcal{B}=W^{k,p}(\mathbb{R})$, with $k\geq3$ and
$1\leq p<\infty$, and suppose that the assumptions of
Theorem~\ref{prop-boundedness} hold. Fix $\sigma>0$, and let
$(\tilde v,\tilde u)$ be the solution of the constant-diffusion
system \eqref{fhn-knwon} with diffusion coefficient $\sigma$
and initial datum
\[
(\tilde v_0,\tilde u_0)\in\mathcal{B}\times\mathcal{B}.
\]
Then, for every $\mu>0$, there exist
$\varepsilon_1,\varepsilon_2$ such that, if
\begin{equation}
\sup_{x\in\mathbb{R}}|\delta(x)-\sigma|<\varepsilon_1,
\end{equation}
and the initial data
$(v_0,u_0)\in\mathcal{B}\times\mathcal{B}$ satisfies
\begin{equation}
\sup_{x\in\mathbb{R}}
\left|
(v_0(x),u_0(x))-(\tilde v_0(x),\tilde u_0(x))
\right|<\varepsilon_2,
\end{equation}
the corresponding solution $(v,u)$ of the variable-diffusion
system \eqref{delta-fhn} satisfies
\begin{equation}
\begin{aligned}
|v(x,t)-\tilde v(x,t)|&\leq\mu,\\
|u(x,t)-\tilde u(x,t)|&\leq\mu,
\end{aligned}
\qquad \forall x\in\mathbb{R},\ t>0.
\end{equation}
\end{theorem}

The paper is organized as follows. Section 3 provides background on the existence of the Fundamental Solution for the heat equation with variable diffusion, along with specific bounds on the Fundamental Solution that are essential for our analysis. Section 4 establishes key well-posedness results, beginning with the proof of local existence, followed by the proof of global existence. Additionally, we state a stability result for the system 
\eqref{delta-fhn}. Finally, section 5 focuses on the proof of the main result, stating first a comparison system that will help us to compare the solutions of the systems \eqref{delta-fhn} and \eqref{fhn-knwon}, analyzing first the linear system and extending to the nonlinear case. 
\section{Fundamental solutions}\label{section-fundametal-solution}

In this section, we introduce the fundamental solution associated
with a heat equation with variable diffusion and state the
properties required in the subsequent analysis. We refer to
\cite{friedman2008partial,pascucci2018parametrix} for the general
theory and proofs. Fix $T>0$ and consider
\begin{equation}\label{general-heat-equation}
\mathbb{L}u
:=\partial_tu-\partial_x\bigl(\tilde a(t,x)\partial_xu\bigr)=0,
\qquad (t,x)\in(0,T]\times\mathbb{R}.
\end{equation}
We impose the following conditions:
\begin{enumerate}
\item[$B_1$]
There is a constant $r\geq1$ such that
\[
r^{-1}\leq\tilde a(t,x)\leq r,
\qquad (t,x)\in[0,T]\times\mathbb{R}.
\]

\item[$B_2$]
The coefficient $\tilde a(t,x)$ together with its spatial
derivative $\partial_x\tilde a(t,x)$ are H\"older continuous in $x$ with
some exponent $\theta\in(0,1)$, uniformly with respect to
$t\in[0,T]$.
\end{enumerate}

\begin{theorem}[Existence of the fundamental solution and
representation formula]\label{theorem-for-generalsolution}
Suppose that $B_1$ and $B_2$ hold. Then $\mathbb{L}$ admits
a fundamental solution $G(t,x;\tau,\xi)$ for
$0\leq\tau<t\leq T$ and $x,\xi\in\mathbb{R}$.

Let $u_0\in BUC(\mathbb{R})$, and let $f$ be bounded and
continuous on $[0,T]\times\mathbb{R}$, locally H\"older
continuous in $x$ uniformly with respect to $t$.
Then the Cauchy problem
\begin{equation}\label{cauchy-general-problem}
\begin{cases}
\mathbb{L}u(t,x)=f(t,x),
&(t,x)\in(0,T]\times\mathbb{R},\\
u(0,x)=u_0(x),
&x\in\mathbb{R},
\end{cases}
\end{equation}
has a unique bounded classical solution for $t>0$,
continuous up to $t=0$, given by
\begin{equation}\label{solution-formula-for-general-heat-equation}
\begin{aligned}
u(t,x)
&=\int_{\mathbb{R}}
G(t,x;0,\xi)u_0(\xi)\,d\xi+\int_0^t\int_{\mathbb{R}}
G(t,x;s,\xi)f(s,\xi)\,d\xi\,ds.
\end{aligned}
\end{equation}
Moreover,
\[
\|u(t,\cdot)\|_{L^\infty(\mathbb{R})}
\leq\|u_0\|_{L^\infty(\mathbb{R})}
+\int_0^t\|f(s,\cdot)\|_{L^\infty(\mathbb{R})}\,ds.
\]
\end{theorem}

The fundamental solution has the following properties for
$0\leq\tau<t\leq T$:

\begin{enumerate}[
    label=\textnormal{Property \arabic*.},
    leftmargin=2.5cm,
    labelsep=0.6em,
    itemsep=1em
]
\item \label{property 1}
Preservation of constants:
\[
\int_{\mathbb{R}}G(t,x;\tau,\xi)\,d\xi=1,
\qquad x\in\mathbb{R}.
\]

\item Composition:
\[
G(t,x;\tau,\xi)
=\int_{\mathbb{R}}
G(t,x;\beta,y)G(\beta,y;\tau,\xi)\,dy,
\qquad \tau<\beta<t.
\]

\item Positivity:
\[
G(t,x;\tau,\xi)>0,
\qquad x,\xi\in\mathbb{R}.
\]
\end{enumerate}

The variables $(\tau,\xi)$ specify the time and position
of the source, while $(t,x)$ specify the time and position
at which its effect is evaluated. When $\tilde a$ is
independent of time, the fundamental solution satisfies
\[
G(t,x;\tau,\xi)=G(t-\tau,x;0,\xi).
\]
In this case, Property~2 gives the semigroup identity for
the associated integral operators. For time-dependent
coefficients, it gives the corresponding composition
identity for the two-parameter evolution operators.

\begin{lemma}[Gaussian bounds and spatial derivatives]
\label{bounds-for-fundamental-solution}
Suppose that $B_1$ and $B_2$ hold. There exist constants
$C,A>0$, depending on $T$, the ellipticity bounds, and
the coefficient bounds and H\"older constants, such that,
for all $0\leq\tau<t\leq T$ and $x,\xi\in\mathbb{R}$,
\[
G(t,x;\tau,\xi)
\leq\frac{C}{(t-\tau)^{1/2}}
\exp\left(-\frac{|x-\xi|^2}{4A(t-\tau)}\right),
\]
and
\[
|\partial_xG(t,x;\tau,\xi)|
\leq\frac{C}{t-\tau}
\exp\left(-\frac{|x-\xi|^2}{4A(t-\tau)}\right),
\]
\end{lemma}

We omit the proof of
Lemma~\ref{bounds-for-fundamental-solution}; these estimates
follow from \cite[Theorem~4.7, estimates
(4.9)--(4.11)]{pascucci2018parametrix}.

\section{Well Posedness}
In this section, we extend the results of \cite{RAUCH197812} by considering variable diffusion in place of the constant diffusion used in the original analysis. Although the introduction of a variable diffusion coefficient introduces additional complexity, the fundamental behavior of the solution remains closely aligned with that of the constant diffusion case. This similarity allows us to apply analogous analytical techniques.
\subsection{Functional Framework}
We consider the nonlinear system
\begin{equation}\label{general-equation-fhn}
    U_t=\partial_x\bigl(A(x)U_x\bigr)+F(U),
    \qquad x\in\mathbb{R},\quad t>0,
\end{equation}
with initial condition $U(0,\cdot)=U_0$. Here
$U=(u_1,\ldots,u_n)$,
$A(x)=\operatorname{diag}(a_1(x),\ldots,a_n(x))$, and
$F:\mathbb{R}^n\to\mathbb{R}^n$ is smooth with $F(0)=0$.
For each $i$, assume that either $a_i\equiv0$ or $a_i$ and
its spatial derivative $a_i'$ are bounded and H\"older
continuous on $\mathbb{R}$ with exponent $\theta\in(0,1)$,
and there exists a constant $r_i\geq1$ such that
\[
r_i^{-1}\leq a_i(x)\leq r_i,
\qquad x\in\mathbb{R}.
\]

For each component with diffusion, let
$G_i(t,x;\tau,\xi)$ be the fundamental solution associated
with $\partial_t-\partial_x(a_i(x)\partial_x)$, as described
in Section~\ref{section-fundametal-solution}. Define
\[
(S_i(t)w)(x)
=\int_{\mathbb{R}}G_i(t,x;0,\xi)w(\xi)\,d\xi,
\qquad t>0,
\]
and set $S_i(0)=I$. For components without diffusion,
set $S_i(t)=I$ for all $t\geq0$. Finally, write
\[
\mathcal{S}(t)
=\operatorname{diag}(S_1(t),\ldots,S_n(t)).
\]
Since the coefficients are independent of time, the
composition property of the fundamental solutions gives
the semigroup identity
\[
\mathcal{S}(t+s)=\mathcal{S}(t)\mathcal{S}(s),
\qquad t,s\geq0.
\]

Let $\mathcal{B}$ be a Banach space of real-valued functions
on $\mathbb{R}$. We equip $\mathcal{B}^n$ with the norm
\[
\|U\|_{\mathcal{B}^n}
=\max_{1\leq i\leq n}\|u_i\|_{\mathcal{B}},
\]
and, for $T>0$, define
\[
\|U\|_{C([0,T];\mathcal{B}^n)}
=\sup_{0\leq t\leq T}\|U(t)\|_{\mathcal{B}^n}.
\]

A mild solution of \eqref{general-equation-fhn} on $[0,T]$
is a function $U\in C([0,T];\mathcal{B}^n)$ satisfying
\begin{equation}\label{solution-formula-with-n-components}
U(t)=\mathcal{S}(t)U_0
+\int_0^t\mathcal{S}(t-s)F(U(s))\,ds.
\end{equation}
For each component with diffusion, this identity takes
the form
\[
\begin{aligned}
u_i(t,x)
&=\int_{\mathbb{R}}
G_i(t,x;0,\xi)u_{0,i}(\xi)\,d\xi+\int_0^t\int_{\mathbb{R}}
G_i(t,x;s,\xi)F_i(U(s,\xi))\,d\xi\,ds.
\end{aligned}
\]
For each component without diffusion, it becomes
\[
u_i(t,x)
=u_{0,i}(x)+\int_0^tF_i(U(s,x))\,ds.
\]
Also, let $K\geq1$ satisfy
\[
\|\mathcal S(t)U\|_{\mathcal B^n}
\leq K\|U\|_{\mathcal B^n},
\qquad U\in\mathcal B^n,\quad 0\leq t\leq1.
\]
\subsection{Local Solvability}

\begin{theorem}[Local existence and uniqueness]
\label{theorem-local-solution}
Suppose that \ref{A.1}, \ref{A.2}, and \ref{A.4} hold.
Then, for every $U_0\in\mathcal{B}^n$, there exists
$t_0\in(0,1]$, depending only on $K$,
$\|U_0\|_{\mathcal{B}^n}$, and the Lipschitz bounds of $F$,
such that \eqref{general-equation-fhn} with initial data
$U_0$ admits a unique mild solution
\[
U\in C([0,t_0];\mathcal{B}^n).
\]
Moreover,
\[
\|U\|_{C([0,t_0];\mathcal{B}^n)}
\leq2K\|U_0\|_{\mathcal{B}^n}.
\]
If $\mathcal{S}(t)$ is contractive on $\mathcal{B}^n$,
one may take $K=1$.
\end{theorem}
\begin{proof}
Set $R(t,\cdot)=\mathcal{S}(t)g$. The fundamental-solution
representation gives the integral equation
\[
U(t)=R(t,\cdot)
+\int_0^t\mathcal{S}(t-s)F(U(s))\,ds.
\]
By \ref{A.1}, evaluation at each spatial point is continuous,
so this identity also has a pointwise interpretation. In
particular, for each diffusive component,
\[
\begin{aligned}
u_i(t,x)
&=\int_{\mathbb{R}}G_i(t,x;0,\xi)g_i(\xi)\,d\xi+\int_0^t\int_{\mathbb{R}}
G_i(t,x;s,\xi)F_i(U(s,\xi))\,d\xi\,ds.
\end{aligned}
\]
For a component without diffusion, the corresponding operator
$S_i(t)$ is the identity.

Suppose first that $\|g\|_{\mathcal{B}\times\mathcal{B}}>0$.
For $0<t_0\leq1$, let
\[
\Omega=
\left\{
U\in C([0,t_0];\mathcal{B}\times\mathcal{B}):
\|U-R\|_{C([0,t_0];\mathcal{B}\times\mathcal{B})}
\leq K\|g\|_{\mathcal{B}\times\mathcal{B}}
\right\}.
\]
By \ref{A.2},
\[
R\in C([0,t_0];\mathcal{B}\times\mathcal{B}),
\qquad
\|R\|_{C([0,t_0];\mathcal{B}\times\mathcal{B})}
\leq K\|g\|_{\mathcal{B}\times\mathcal{B}}.
\]
Thus, $\Omega$ is a nonempty closed subset of the Banach space
$C([0,t_0];\mathcal{B}\times\mathcal{B})$. Moreover, for
$U\in\Omega$,
\[
\begin{aligned}
\|U\|_{C([0,t_0];\mathcal{B}\times\mathcal{B})}
&\leq
\|U-R\|_{C([0,t_0];\mathcal{B}\times\mathcal{B})}
+\|R\|_{C([0,t_0];\mathcal{B}\times\mathcal{B})}\\
&\leq2K\|g\|_{\mathcal{B}\times\mathcal{B}}.
\end{aligned}
\]
By \ref{A.4}, there exists a constant $k\geq0$, independent
of $t_0\leq1$, such that
\begin{equation}\label{lipschitz-property-middleproof-localsolvae}
\|F(V)-F(W)\|_{\mathcal{B}\times\mathcal{B}}
\leq k\|V-W\|_{\mathcal{B}\times\mathcal{B}}
\end{equation}
whenever
$\|V\|_{\mathcal{B}\times\mathcal{B}},
\|W\|_{\mathcal{B}\times\mathcal{B}}
\leq2K\|g\|_{\mathcal{B}\times\mathcal{B}}$.
Since $F(0)=0$, this also gives
\[
\|F(V)\|_{\mathcal{B}\times\mathcal{B}}
\leq k\|V\|_{\mathcal{B}\times\mathcal{B}}
\]
on the same ball.

Choose $t_0\in(0,1]$ such that $Kkt_0\leq1/4$, and define
\[
\Gamma U(t)=R(t,\cdot)
+\int_0^t\mathcal{S}(t-s)F(U(s))\,ds.
\]
We first verify that $\Gamma U$ belongs to
$C([0,t_0];\mathcal{B}\times\mathcal{B})$ for $U\in\Omega$.
By \ref{A.4}, the map $s\mapsto F(U(s))$ is continuous.
The strong continuity and boundedness in \ref{A.2} imply that
\[
(t,s)\longmapsto\mathcal{S}(t-s)F(U(s))
\]
is continuous on 
$0\leq s\leq t\leq t_0$. Consequently,
\[
t\longmapsto
\int_0^t\mathcal{S}(t-s)F(U(s))\,ds
\]
is continuous. Together with the continuity of $R$, this proves
that $\Gamma U\in C([0,t_0];\mathcal{B}\times\mathcal{B})$.

We next show that $\Gamma$ maps $\Omega$ into itself. By
\ref{A.2} and \eqref{lipschitz-property-middleproof-localsolvae},
for $U\in\Omega$ and $0\leq t\leq t_0$,
\[
\begin{aligned}
\|\Gamma U(t)-R(t,\cdot)\|_{\mathcal{B}\times\mathcal{B}}
&\leq K\int_0^t
\|F(U(s))\|_{\mathcal{B}\times\mathcal{B}}\,ds\\
&\leq Kk\int_0^t
\|U(s)\|_{\mathcal{B}\times\mathcal{B}}\,ds\\
&\leq2K^2kt_0\|g\|_{\mathcal{B}\times\mathcal{B}}\\
&\leq\frac{K}{2}\|g\|_{\mathcal{B}\times\mathcal{B}}.
\end{aligned}
\]
Taking the supremum over $t\in[0,t_0]$, we obtain
$\Gamma(\Omega)\subseteq\Omega$.

For $U,V\in\Omega$, the same assumptions give
\[
\begin{aligned}
\|\Gamma U(t)-\Gamma V(t)\|_{\mathcal{B}\times\mathcal{B}}
&\leq K\int_0^t
\|F(U(s))-F(V(s))\|_{\mathcal{B}\times\mathcal{B}}\,ds\\
&\leq Kk\int_0^t
\|U(s)-V(s)\|_{\mathcal{B}\times\mathcal{B}}\,ds\\
&\leq Kkt_0
\|U-V\|_{C([0,t_0];\mathcal{B}\times\mathcal{B})}.
\end{aligned}
\]
Therefore,
\[
\|\Gamma U-\Gamma V\|_{C([0,t_0];\mathcal{B}\times\mathcal{B})}
\leq\frac14
\|U-V\|_{C([0,t_0];\mathcal{B}\times\mathcal{B})}.
\]
Thus, $\Gamma$ is a contraction on the complete set $\Omega$.
Banach's fixed-point theorem yields a unique fixed point
$U\in\Omega$, which is a mild solution of
\eqref{general-equation-fhn}. Since $\mathcal{S}(0)=I$,
the integral equation gives $U(0)=g$. Furthermore,
\[
\|U\|_{C([0,t_0];\mathcal{B}\times\mathcal{B})}
\leq2K\|g\|_{\mathcal{B}\times\mathcal{B}}.
\]

To establish uniqueness among all mild solutions in
$C([0,t_0];\mathcal{B}\times\mathcal{B})$, let $U$ and $V$
be two such solutions with the same initial data. Their ranges
are contained in a common bounded ball of
$\mathcal{B}\times\mathcal{B}$. By \ref{A.4}, there exists a
Lipschitz constant $k_*$ on that ball. Using \ref{A.2}, we obtain
\[
\|U(t)-V(t)\|_{\mathcal{B}\times\mathcal{B}}
\leq Kk_*\int_0^t
\|U(s)-V(s)\|_{\mathcal{B}\times\mathcal{B}}\,ds.
\]
Gronwall's inequality implies $U=V$.

Finally, if $g=0$, the function $U=0$ is a mild solution because
$F(0)=0$, and the preceding uniqueness argument applies. This
completes the proof.
\end{proof}
To demonstrate the global solvability of the solution, we use the technique of contracting rectangles presented initially on \cite{RAUCH197812}. The key property of contracting rectangles is that they enable the definition of nonlinear functionals, which decrease over time for the solutions. Let $|\cdot|_R$ be the norm on $\mathbb{R}^2$ defined by such a rectangle in the usual way:
$$
|U|_R=\inf \{t \geqslant 0: U \in t R\} .
$$
And we define a continuous function $\mathscr{V}_R\colon BC^0(\mathbb{R})\times BC^0(\mathbb{R})
\to\mathbb{R}$ by
\begin{equation}\label{lyapunov-functional}
    \mathscr{V}_R(W)=\sup _{x \in \R}|W(x)|_R .
\end{equation}
Also for a function $\Psi: \mathbb{R} \rightarrow \mathbb{R}$, the upper Dini derivate by
\begin{equation}\label{dini-derivative}
 \bar{D} \Psi(T)=\limsup _{h \rightarrow 0} \frac{\Psi(T+h)-\Psi(T)}{h}.   
\end{equation}
We now present the following lemma, which demonstrates how the functional $\mathscr{V}_R$ decreases over time. The proof closely follows the one given in \cite{RAUCH197812}, but we include it here to highlight that the result remains valid even with variable diffusion.
\begin{lemma}\label{basic-lemma}
Let $F(U)$ be a vector field on $\mathbb{R}^n$, and let $R$ be a rectangle with $0 \in \operatorname{int}(R)$. Suppose that $U \in C\left((T-\eta, T+\eta), C_0(\mathbb{R})\right)$ is a smooth solution of \eqref{delta-fhn} for $|t-T|<\eta$ and that $\mathscr{V}_R(U(T))=s$. If there is an $\eta>0$ such that for any $W \in \partial(s R)$ and $n(W)$ normal to $\partial(s R)$ at $W$ we have $F(W) \cdot n(W)<-\eta$, then
\begin{equation}
\bar{D} \mathscr{V}_R(U(T)) \leqslant-(2 \eta / s L) \mathscr{V}_R(U(T))
\end{equation}
where $L$ is the length of the shortest side of $R$.
\end{lemma}
\begin{proof}
We focus on the case where the solution \( U \) reaches the boundary, because if this does not occur, since \( U(t, \cdot) \in C_0(\R) \), the solution will remain in the interior of \( R \) for \( t \in (T - \eta, T + \eta) \). Additionally, without loss of generality, we assume $\mathscr{V}_R(U(T)) = 1$ and we will carry out the proof for \( n = 2 \). Let $U=(v,u), F= (f_1,f_2)$ and let $R$ be defined by the inequalities $-r_1 \leq v \leq r_2$ and $-s_1\leq u \leq s_2$. Let's do the calculations for $v$, for $u$ is analogue. So we have that $-r_1 \leq v(T,x) \leq r_2$ for all $x$. We say that $U(t,x)$ is in the right hand face if $v(t,x) = r_2$. So, $v(T,x) \leq r_2$ for all $x$ with equality at $x=\bar{x}$, this means, $v_x(T,\bar{x}) = 0$ and since $\delta(x) > 0$ for all $x$, $\delta v_{xx}(T,\bar{x}) \leq 0$. Thus
$$
v_t = (\delta v_x)_x + f_1(v)  < 0
$$
so that $v(T+h,\bar{x}) < r_2 - \eta h$ for small $h$. By continuity, this holds for all $x$ in a neighborhood of $\bar{x}$. For the left side of $R$ we get $v(T+h,\bar{x}) < r_1 - \eta h$, and for the top and bot side we get a similar inequality. So taking in account all this inequalities and since the $\partial R$ and the set of outward pointing unit normals to $R$ at the boundary is compact, there exists a neighborhood $N(\bar{x})$ such that for $y \in N(\bar{x})$ we have
$$
U(T+h, y) \in\left(1-h \eta\left(\min _{1 \leqslant j \leqslant n}\left(r_j, s_j\right)\right)^{-1}\right) R.
$$
Thus, for sufficiently small $h$ we have $U(T+h, x) \in\left(1-h \eta\left(\min _{1 \leqslant j \leqslant n}\left(r_j, s_j\right)\right)^{-1}\right) R$ for every $x\in \R$. Which is equivalent that the functional \eqref{lyapunov-functional} satisfy 
$$
\mathscr{V}_R(U(T+h)) \leqslant 1-(2 h \eta / L),
$$
where $L = \min _{1 \leqslant j \leqslant n}\left(r_j, s_j\right)$. Then, according to \eqref{dini-derivative} we get
\begin{equation}
\bar{D} \mathscr{V}_R(U(T)) \leqslant-(2 \eta / s L) \mathscr{V}_R(U(T)),
\end{equation}
which completes the proof.
\end{proof}
\subsection{Proof of Theorem \ref{prop-boundedness}}
Lemma \ref{basic-lemma} allow us to prove the first main result of this work.
\begin{proof}
First we prove the global solvability and after we focus on the stability estimate. So, for this, assume that there exists a contracting rectangle $R$ such that $\mathscr{V}_R(U_0)<1$ and consider $U \in C\left(\left[0, t_0\right] ; \mathcal{B} \times \mathcal{B}\right)$, which is the local solution provided by Theorem \ref{theorem-local-solution} , satisfying the initial condition $U(0)=U_0$. 

Now, since the function $\mathscr{V}_R(U(t))$, defined for $t \in\left[0, t_0\right]$, is continuous and satisfies $\mathscr{V}_R(U_0)<1$, along with the fact that its upper right Dini derivative $\bar{D} \mathscr{V}_R(U(t))$ remains negative whenever $\mathscr{V}_R(t)$ reaches 1 , we conclude that $\mathscr{V}_R(U(t))$ stays strictly below 1 for all $t \in\left[0, t_0\right]$.
Specifically, evaluating at $t=t_0 / 2$, we have $\mathscr{V}_R\left(U(t_0 / 2\right))<1$. This observation allows us to restart the argument, extending the solution $U$ over the larger interval $\left[0, t_0+t_0 / 2\right]$ while maintaining $\mathscr{V}_R(U(t))<1$ throughout.
Repeating this procedure indefinitely in increments of $t_0 / 2$, we conclude that the solution $U$ remains well defined for all $t \geq 0$. Proving the first statement of the Theorem.

%We now present a key result that provides an explicit estimate on the smallness condition required for the initial data $U^0 = (v_0, u_0)$ to ensure stability. This result will serve as a crucial ingredient in the subsequent sections.
%\begin{theorem}\label{stability-deltafhn}
%    Suppose $a>1 / \gamma$ and $U^0 \in C_0(\mathbb{R})$. Then given $m>0$ there exists a rectangle $R$  such that $\mathscr{V}_R(U_0) < m$. Moreover, there are nonnegative constants $c, K$ such that $\|U(t)\|_{\infty} \leqslant K e^{-c t}$, for all $t \geqslant 0$.
%\end{theorem}
%\begin{proof}
We move on to establish the stability estimate, for this, we need to show the existence of a contracting rectangle $R$. We verify that the vector field $H(U, x, t)$ satisfies Definition \ref{definition-contracting-rectangle}. The vector field is given by

$$
H(U, x, t) = \big(f(v) - u + \delta^{\prime} v_x, \epsilon(v - \gamma u) \big).
$$

Since the definition applies to every $U \in \partial R$, we have $v_x = 0$ on the boundary. Consequently, the vector field simplifies to

$$
H(U, x, t) = \big(f(v) - u, \epsilon(v - \gamma u) \big).
$$

Thus, since $a> 1/\gamma$ by [Lemma 3.4, \cite{RAUCH197812}], for any given $m>0$ there exists a contracting rectangle $R$ under the flow of the vector field $H$. Next, since $U_0 \in C_0(\mathbb{R})$, we can select $m \ll 1$ such that $\mathscr{V}_R(U_0) < m.$ Applying Lemma \ref{basic-lemma}, we obtain the differential inequality:

$$
\bar{D} \mathscr{V}_R(U(T)) \leqslant -\frac{2 \eta}{m L} \mathscr{V}_R(U(T)).
$$

Solving this inequality gives

$$
\mathscr{V}_R(U(t)) \leq e^{-\frac{2 \eta t}{m L}}.
$$

This result implies the existence of nonnegative constants $c, K$ such that

$$
\| U(t) \|_{\infty} \leqslant K e^{-c t}, \quad \forall t \geq 0.
$$

Thus, the contraction property ensures the exponential decay of $\| U(t) \|_{\infty}$, completing the proof.

\end{proof}
We end this section by showing uniform bounds on the solution.
\begin{proposition}[Properties of the solution]
\label{proposition-properties-solution}
Let $\mathcal{B}=W^{k,p}(\mathbb{R})$, with $k\geq 3$ and
$1\leq p<\infty$, and let $U^0=(v_0,u_0)\in\mathcal{B}\times\mathcal{B}$.
Under the assumptions of Theorem \ref{prop-boundedness}, let
\[
U=(v,u)\in C\bigl([0,\infty);
\mathcal{B}\times\mathcal{B}\bigr)
\]
be the unique global solution of \eqref{delta-fhn}.
Then there exist constants $C_1,C_2>0$ such that
\[
|v(x,t)|\leq C_1,
\qquad
|\partial_x v(x,t)|\leq C_2,
\qquad x\in\mathbb{R},\quad t\geq 0.
\]
\end{proposition}

\begin{proof}
Throughout the proof, $C$ denotes a positive constant that may change
from line to line but is independent of $t$.
By Theorem \ref{prop-boundedness},
\[
\sup_{t\geq 0}
\bigl(\|v(t)\|_{L^\infty(\mathbb{R})}
+\|u(t)\|_{L^\infty(\mathbb{R})}\bigr)\leq C.
\]
This proves the first estimate and, since $f$ is a polynomial, also gives
\[
\sup_{t\geq 0}
\|f(v(t))-u(t)\|_{L^\infty(\mathbb{R})}\leq C.
\]

On $[0,1]$, continuity in $W^{k,p}(\mathbb{R})$ and the Sobolev
embedding yield
\[
\sup_{0\leq t\leq 1}\|\partial_x v(t)\|_{L^\infty(\mathbb{R})}
\leq C\sup_{0\leq t\leq 1}\|v(t)\|_{W^{k,p}(\mathbb{R})}
<\infty.
\]

For $t\geq 1$, we use the representation formula
\eqref{solution-formula-for-general-heat-equation} on $[t-1,t]$:
\[
\begin{aligned}
v(t,x)
={}&\int_{\mathbb{R}}G(t,x;t-1,\xi)v(t-1,\xi)\,d\xi\\
&+\int_{t-1}^{t}\int_{\mathbb{R}}
G(t,x;s,\xi)\bigl(f(v(s,\xi))-u(s,\xi)\bigr)\,d\xi\,ds.
\end{aligned}
\]
Since $\delta$ is independent of time, Lemma
\ref{bounds-for-fundamental-solution}, applied on intervals of length
one, gives
\[
\sup_{x\in\mathbb{R}}
\int_{\mathbb{R}}|\partial_xG(t,x;s,\xi)|\,d\xi
\leq C(t-s)^{-1/2},
\qquad 0<t-s\leq 1,
\]
with $C$ independent of the initial time.
The integrability of this bound allows us to differentiate the
representation formula once in $x$. Thus,
\[
\begin{aligned}
\|\partial_x v(t)\|_{L^\infty(\mathbb{R})}
&\leq C\|v(t-1)\|_{L^\infty(\mathbb{R})}\\
&\quad+C\int_{t-1}^{t}(t-s)^{-1/2}
\|f(v(s))-u(s)\|_{L^\infty(\mathbb{R})}\,ds\\
&\leq C\left(1+\int_0^1r^{-1/2}\,dr\right)
\leq C.
\end{aligned}
\]
Combining the bounds on $[0,1]$ and $[1,\infty)$ proves the second
estimate.
\end{proof}
The constant-diffusion FitzHugh--Nagumo system also satisfies
properties analogous to those stated in Proposition
\ref{proposition-properties-solution}. In particular, its solution
$\widetilde v$ and its first spatial derivative
$\partial_x\widetilde v$ are uniformly bounded in space and time;
see \cite{cerpa2023approximation}.
\section{Continuity Result}
The objective of this section is to establish Theorem \ref{thm-approximation} by analyzing how the solution of \eqref{delta-fhn} compare to the solution of \eqref{fhn-knwon} given that
\begin{equation}\label{condition-initial-diffusion}
\sup_{x\in \R}|\delta(x) - \sigma| < \epsilon_1, \quad\sup_{x\in \R}|(v_0(x),u_0(x))-(\tilde{v}_0(x), \tilde{u}_0(x))| < \epsilon_3,
\end{equation}
for $\epsilon_1, \epsilon_3>0$ chosen small enough. We first focus on the linearized version of the system to derive estimates that depend on its solution and impose certain constraints on the parameters. Using these preliminary results, we then extend the analysis to the full nonlinear case and show that, under appropriate assumptions, the difference between the two solutions remains uniformly controlled for all time. For this we do the natural change of variables 
\begin{equation}
    E_v =  v - \tilde{v} \quad \text{and} \quad E_u = u - \tilde{u},
\end{equation}
therefore, we get the comparison system,
\begin{equation}\label{difference-system}
\left\{\begin{aligned}
\partial_t E_v = \sigma \partial_x^2 E_v &+ (\varphi_1 - a) E_v + \varphi_2 E_v^2 - E_v^3 - \varphi_3 - E_u, \\
\partial_t E_u  =\rho \partial_x^2 E_u &+ \varepsilon \left(E_v-\gamma E_u\right) , \\
E_v(x, 0)=E_v(0), & E_u(x, 0)=E_u(0),
\end{aligned}\right.
\end{equation}
where
\begin{align}
    \varphi_1 &= 2(1+a)\tilde{v} - 3\tilde{v}^2 ,\\
    \varphi_2 &= 1+a-3\tilde{v}\\
    \varphi_3 &= \left((\delta - \sigma)v_x \right)_x.
\end{align}
\subsection{\texorpdfstring{Linear estimate of the system \eqref{difference-system}}{Linear estimate of the system ()}}
We prove an approximation result for the linear problem
\begin{equation}\label{linear-error-equation}
\left\{\begin{aligned}
\partial_t F_v-\sigma\partial_x^2 F_v & =(\varphi_1-a) F_v+\varphi_3-F_u, \\
\partial_t F_u-\rho \partial_x^2 F_u &= \varepsilon \left(F_v-\gamma F_u\right)  , \\
F_v(x, 0)=E_v(0), & F_u(x, 0)=E_u(0),
\end{aligned}\right.
\end{equation}
\begin{lemma}[Linear comparison system]
\label{lemma-linear-error}
Let $\gamma>0$, $\sigma>0$, and $a\in(0,\tfrac12)$ satisfy
$a>1/\gamma$. Let $\varepsilon_1,\mu_1,\mu_2>0$, and suppose that
\[
\|\delta-\sigma\|_{L^\infty(\mathbb R)}<\varepsilon_1.
\]
Assume that the solution $(\tilde v,\tilde u)$
of \eqref{fhn-knwon} satisfies
\begin{equation}\label{condition-for-contraction}
\alpha(T)
:=2(1+a)\|\tilde v\|_Y+3\|\tilde v\|_Y^2
<a-\frac{1}{\gamma}.
\end{equation}
Suppose further that the initial data satisfy
\[
\|F_v(0,\cdot)\|_{L^\infty(\mathbb R)}\leq\mu_1,
\qquad
\|F_u(0,\cdot)\|_{L^\infty(\mathbb R)}\leq\mu_2.
\]

Then the initial value problem \eqref{linear-error-equation}
admits a solution
\[
(F_v,F_u)\in C([0,T],\mathcal B\times\mathcal B),
\]
and there exists a constant $C>0$, independent of
$\mu_1,\mu_2,\varepsilon_1$, such that
\begin{equation}\label{bounds-linear-equation}
\begin{aligned}
|F_v(t,x)|&\leq C(\mu_1+\mu_2+\varepsilon_1),\\
|F_u(t,x)|&\leq C(\mu_1+\mu_2+\varepsilon_1),
\end{aligned}
\qquad
\forall x\in\mathbb R,\quad 0\leq t\leq T.
\end{equation}
\end{lemma}

Since \eqref{linear-error-equation} has constant diffusion
coefficient $\sigma>0$, the following proof uses the heat kernel
\[
G_\sigma(t,x)=\frac{1}{\sqrt{4\pi\sigma t}}
\exp\left(-\frac{x^2}{4\sigma t}\right),
\qquad t>0,
\]
\begin{proof}
    First, we multiply the first equation on \eqref{linear-error-equation} by $e^{a t}$ and we get
\begin{equation*}
\partial_t e^{a t} F_v-\sigma\partial_x^2 e^{a t} F_u=\varphi_1 e^{a t} F_v+\varphi_3e^{a t} -  e^{a t} F_u.
\end{equation*}   
And we multiply the second equation on \eqref{linear-error-equation} by $e^{\varepsilon \gamma t}$ and we get
$$
\partial_t e^{\varepsilon \gamma t} F_u-\rho\partial_x^2 e^{\varepsilon \gamma t} F_u=\varepsilon  e^{\varepsilon \gamma t} F_v
$$
And by virtue of Duhamel's principle, we get
$$
F_v(t,x)=G_{\sigma}(t) * F_v(0,\cdot)+\int_0^t e^{-a(t-\tau)} G_{\sigma}(t-\tau) *\left(\varphi_1 F_v+\varphi_3 - F_u\right) d \tau .
$$
For the second equation, we obtain
$$
F_u(t,x)=G_{\rho}(t) * F_u(0,\cdot)+\varepsilon \int_0^t e^{-\varepsilon \gamma(t-\tau)} G_{\rho}(t-\tau) *F_v d \tau .
$$
Let us consider the following iterative procedure. Let $\mu_1, \mu_2$ small and set $F_v^{(0)}=\mu_1, F_u^{(0)}=\mu_2$ and define
$$
\begin{aligned}
& F_v^{(k+1)}(t,x)=G_{\sigma}(t) * F_v(0,\cdot)+\int_0^t e^{-a(t-\tau)} G_{\sigma}(t-\tau) *\left(\varphi_1 F_v^{(k)}+\varphi_3 - F_u^{(k)}\right) d \tau \\
& F_u^{(k+1)}(t,x)=G_{\rho}(t) * F_u(0,\cdot)+\varepsilon \int_0^t e^{-\varepsilon \gamma(t-\tau)} G_{\rho}(t-\tau) *F_v^{(k)} d \tau,
\end{aligned}
$$
where $\left(F_v^{(k)}, F_u^{(k)}\right) \in B C^0 \times B C^0$ imply that $\left(F_v^{(k+1)}, F_u^{(k+1)}\right) \in B C^0 \times B C^0$. The next step is to look at the convergence of the sequences $\left\{F_v^{(k)}\right\},\left\{F_u^{(k)}\right\}$ in the space $C\left([0, T] ; B C^0\right)$. We consider the norm in $C\left([0, T] ; B C^0\right)$ given by (25) and estimate the difference between two consecutive terms. For $\left\{F_v^{(k)}\right\}$, it holds

$$
\begin{aligned}
\left\|F_v^{(k+1)}-F_v^{(k)}\right\|_Y &\leq  \sup _{0<t<T} \int_0^t e^{-a(t-\tau)} \left\|G_{\sigma}(t-\tau) *\left|\varphi_1(\tau,\cdot)\right|\right\|_{L_x^{\infty}} d \tau\left\|F_v^{(k)}-F_v^{(k-1)}\right\|_Y \\
&+ \sup _{0<t<T} \int_0^t e^{-a(t-\tau)}\left\|F_u^{(k)}-F_u^{(k-1)}\right\|_Y \\
& \leq \sup _{0<t<T} \int_0^t e^{-a(t-\tau)} \left\|G_{\sigma}(t-\tau) *\left|\varphi_1(\tau,\cdot)\right|\right\|_{L_x^{\infty}} d \tau\left\|F_v^{(k)}-F_v^{(k-1)}\right\|_Y \\
&+ \sup _{0<t<T} \int_0^t e^{-a(t-\tau)}\dfrac{1}{\gamma }\left\|F_v^{(k)}-F_v^{(k-1)}\right\|_Y .
\end{aligned}
$$
Now,
\begin{align*}
\left\|G_{\sigma}(t-\tau) *\left|\varphi_1(\tau,\cdot)\right|\right\|_{L_x^{\infty}} &\leq  \left\| \varphi_1(\tau,\cdot)\right\|_{L_x^{\infty}}\\
&\leq 2(1+a)\norm{\tilde{v}}_Y + 3\norm{\tilde{v}}_Y^2.
\end{align*}
Thus,
\begin{equation*}
\begin{aligned}
\left\|F_v^{(k+1)}-F_v^{(k)}\right\|_Y  &\leq\left\|F_v^{(k)}-F_v^{(k-1)}\right\|_Y\left(2(1+a)\norm{\tilde{v}}_Y + 3\norm{\tilde{v}}_Y^2+\dfrac{1}{\gamma}\right)\sup _{0 \leq t \leq T} \int_0^t e^{-a(t-\tau)} d \tau \\
&\leq \left\|F_v^{(k)}-F_v^{(k-1)}\right\|_Y \dfrac{2(1+a)\norm{\tilde{v}}_Y + 3\norm{\tilde{v}}_Y^2+\dfrac{1}{\gamma}}{a}.
\end{aligned}
\end{equation*}
Defining $\alpha(T)=\dfrac{2(1+a)\norm{\tilde{v}}_Y + 3\norm{\tilde{v}}_Y^2+\dfrac{1}{\gamma}}{a}$ we obtained that
$$
\left\|F_v^{(k+1)}-F_t^{(k)}\right\|_Y \leq \alpha(T)\left\|F_t^{(k)}-F_v^{(k-1)}\right\|_Y .
$$
Analogously, for $\left\{F_u^{(k)}\right\}$, we get
$$
\left\|F_u^{(k+1)}-F_u^{(k)}\right\|_Y \leq \varepsilon \sup _{0<t<T} \int_0^t e^{-\varepsilon \gamma(t-\tau)} d \tau\left\|F_v^{(k+1)}-F_v^{(k)}\right\|_Y \leq \dfrac{1}{\gamma}\left\|F_v^{(k+1)}-F_v^{(k)}\right\|_Y .
$$
Now, in order to state the convergence of the sequences $\left\{F_v^{(k)}\right\},\left\{F_u^{(k)}\right\}$ in $C\left([0, T] ; B C^0\right)$, use that for $m>n$
$$
\begin{gathered}
\left\|F_v^{(m)}-F_v^{(n)}\right\|_Y \leq \sum_{k=n}^{m-1}\left\|F_t^{(k+1)}-F_v^{(k)}\right\|_Y \leq \sum_{k=n}^{m-1} \alpha(T)^k\left\|F_v^{(1)}-F_v^{(0)}\right\|_Y, \\
\left\|F_u^{(m)}-F_u^{(n)}\right\|_Y \leq \sum_{k=n}^{m-1}\left\|F_u^{(k+1)}-F_u^{(k)}\right\|_Y \leq \sum_{k=n}^{m-1} \dfrac{1}{\gamma}\left\|F_v^{(k+1)}-F_v^{(k)}\right\|_Y \leq \sum_{k=n}^{m-1} \dfrac{1}{\gamma} \alpha(T)^k\left\|F_v^{(1)}-F_v^{(0)}\right\|_Y .
\end{gathered}
$$
Condition \eqref{condition-for-contraction} guarantees that we have a contraction mapping hence the sequence $\left\{\left(F_v^{(k)}, F_{u}^{(k)}\right)\right\}_k$ converges in $C\left([0, T] ; B C^0 \times B C^0\right)$ and uniqueness of solution [Theorem 5, \cite{cerpa2023approximation}] imply that $\left(F_v, F_u\right) \in C([0, \infty), \mathcal{B} \times \mathcal{B})$. To conclude the last part of the proof of the Lemma we need to prove the estimate \eqref{bounds-linear-equation} which can be obtained by bounding in a proper manner $\left\|F_v^{(1)}-F_v^{(0)}\right\|_Y$.

\begin{align*}
\left\|F_v^{(1)}-F_v^{(0)}\right\|_Y&=\left\| G_{\sigma}(t) * F_v(0, \cdot) + \Int_0^t e^{-a(t-\tau)} G_{\sigma}(t-\tau) *\left(\varphi_1 F_v^{(0)}+\varphi_3-F_u^{(0)}\right) d \tau  - F_v^{(0)}\right\|_Y \\
&\leq \tilde{C}(\mu_1,\mu_2) + \left\|\int_0^t e^{(-a)(t-\tau)} G_{\sigma}(t-\tau) * \varphi_3 d \tau\right\|_Y\\
&\leq \tilde{C}(\mu_1,\mu_2) + \sup _{0 \leq t \leq T} \int_0^t e^{-a (t-\tau)}\left\| G_{\sigma}(t-\tau) \ast \varphi_3(\tau, \cdot)\right\|_{L^{\infty}_x} d \tau  \\
&\leq \tilde{C}(\mu_1,\mu_2)
+ \sup_{0 \leq t \leq T}
\int_0^t e^{-a(t-\tau)}
\left\|
\partial_x G_{\sigma}(t-\tau)
* \bigl((\delta-\sigma)v_x(\tau,\cdot)\bigr)
\right\|_{L^\infty_x}\,d\tau \\
&\leq \tilde{C}(\mu_1,\mu_2)
+ \sup_{0 \leq t \leq T}
\int_0^t e^{-a(t-\tau)}
\left\|\partial_x G_{\sigma}(t-\tau)\right\|_{L^1_x}
\left\|(\delta-\sigma)v_x(\tau,\cdot)\right\|_{L^\infty_x}
\,d\tau \\
&\leq \tilde{C}(\mu_1,\mu_2)
+ \frac{\|\delta-\sigma\|_{L^\infty_x}\|v_x\|_Y}
{\sqrt{\pi\sigma}}
\sup_{0 \leq t \leq T}
\int_0^t
\frac{e^{-a(t-\tau)}}{\sqrt{t-\tau}}\,d\tau \\
&= \tilde{C}(\mu_1,\mu_2)
+ \frac{\|\delta-\sigma\|_{L^\infty_x}\|v_x\|_Y}
{\sqrt{\pi\sigma}}
\sup_{0 \leq t \leq T}
\int_0^t e^{-as}s^{-1/2}\,ds \\
&\leq \tilde{C}(\mu_1,\mu_2)
+ \frac{\|\delta-\sigma\|_{L^\infty_x}\|v_x\|_Y}
{\sqrt{\pi\sigma}}
\int_0^\infty e^{-as}s^{-1/2}\,ds \\
&= \tilde{C}(\mu_1,\mu_2)
+ \frac{\|\delta-\sigma\|_{L^\infty_x}\|v_x\|_Y}
{\sqrt{a\sigma}} \\
&\leq \tilde{C}(\mu_1,\mu_2)
+ \frac{\varepsilon_1 M}{\sqrt{a\sigma}},
\end{align*}
using Theorem \ref{prop-boundedness}, Proposition \ref{proposition-properties-solution}, and the condition $\delta$ close to $\sigma$ together with the smallness condition on the derivative of $\delta$, we get
\begin{equation}\label{proof-bound-for-sequence}
 \left\|F_v^{(1)}-F_v^{(0)}\right\|_Y \leq \epsilon C. 
\end{equation}
Finally, using the convergence of the sequence $\left\{F_v^{(k)}\right\}$, the estimate \eqref{proof-bound-for-sequence} and that $\alpha=\alpha(T)<1$ we conclude
$$
\begin{aligned}
\left\|F_v-F_v^{(0)}\right\|_Y & \leq\left\|F_v-F_v^{(N+1)}\right\|_Y+\sum_{k=0}^N\left\|F_v^{(k+1)}-F_v^{(k)}\right\|_Y \\
& \leq\left\|F_v-F_v^{(N+1)}\right\|_Y+\sum_{k=0}^N \alpha^k\left\|F_v^{(1)}-F_v^{(0)}\right\|_Y \\
& \leq\left\|F_v-F_v^{(N+1)}\right\|_Y+\frac{1}{1-\alpha} \varepsilon C .
\end{aligned}
$$

Therefore, by taking the limit as $N \rightarrow \infty$ we get the first estimate in \eqref{bounds-linear-equation}. Analogously for $F_u$ we get
$$
\begin{aligned}
\left\|F_u-F_u^{(0)}\right\|_Y & \leq\left\|F_u-F_u^{(N+1)}\right\|_Y+\sum_{k=0}^N\left\|F_u^{(k+1)}-F_u^{(k)}\right\|_Y \\
& \leq\left\|F_u-F_u^{(N+1)}\right\|_Y+\sum_{k=0}^N \frac{1}{\gamma}\left\|F_v^{(k+1)}-F_v^{(k)}\right\|_Y \\
& \leq\left\|F_u-F_u^{(N+1)}\right\|_Y+\frac{1}{\gamma} \frac{\alpha}{1-\alpha} \varepsilon C  .
\end{aligned}
$$

Taking the limit as $N \rightarrow \infty$ we get the second part of \eqref{bounds-linear-equation}.
\end{proof}

\begin{remark}
Since $a>\gamma^{-1}$, Lemma~3.4 in \cite{RAUCH197812} provides a rectangle $R$
containing the origin such that $\tau R$ is contracting for all
sufficiently small $\tau>0$. Let
\[
M_R=\max_{(v,u)\in R}|v|.
\]
If the range of the initial data $(\tilde v_0,\tilde u_0)$ is
contained in $\tau R$, invariance yields
\[
\|\tilde v\|_Y\leq\tau M_R.
\]
Consequently, choosing $\tau>0$ sufficiently small gives
\[
\alpha(T)
\leq2(1+a)\tau M_R+3\tau^2M_R^2
<a-\frac{1}{\gamma},
\]
uniformly for every $T>0$.
\end{remark}
\subsection{\texorpdfstring{Nonlinear estimate of the system \eqref{difference-system}}{Nonlinear estimate of the system ()}}
The next Lemma show the error estimate for the Non-linear equation \eqref{difference-system}. We omit the proof since it is analogous as the proof given in Lemma 15 in \cite{cerpa2023approximation}.
\begin{lemma}\label{lemma-nonlinear-error}
(Non-linear estimate for the comparison system). Let $B =W^{k, p}(\mathbb{R})$ with $k \geq 3, p \geq 1$,$\varepsilon>0$. For some $T>0$ let $(v, u) \in C([0, T] ; \mathcal{B} \times \mathcal{B})$ be a solution of \eqref{delta-fhn}, let $\varphi_1, \varphi_2$ be given by Lemma \ref{lemma-linear-error} and let $\left(F_v, F_u\right) \in C([0, T] ; \mathcal{B} \times \mathcal{B})$ be the corresponding solution of \eqref{linear-error-equation}.
Given $\mu>0$ there exist constants $C_1(\mu), C_2(\mu), C_3(\mu)>0$ such that if
\begin{itemize}
\item $\|v\|_Y \leq C_1(\mu)$,
\item $\left\|\varphi_1\right\|_Y \leq C_2(\mu)$,
\item $\left\|\varphi_2\right\|_Y \leq 3C_1(\mu)+ a+1$,
\item $\left|F_v(t,x)\right| \leq C_3(\mu),\left|F_u(t,x)\right| \leq \mu / 2$, for $x \in \R, 0 \leq t \leq T$,
\end{itemize}
then, there exists a unique solution $\left(E_v, E_u\right) \in C([0, T] ; \mathcal{B} \times \mathcal{B})$ of \eqref{difference-system} that also satisfies the estimate
$$
\left|E_v(t,x)\right| \leq \mu, \quad\left|E_u(t,x)\right| \leq \mu, \quad \forall x \in \R, \quad 0 \leq t \leq T .
$$   
\end{lemma}
\subsection{Proof of Theorem \ref{thm-approximation}}
The existence of a solution $ (v, u) $ to system \eqref{delta-fhn} in $ C([0,\infty), \mathcal{B} \times \mathcal{B}) $ for an initial condition $ (v_0, u_0) \in \mathcal{B} \times \mathcal{B} $, which lies within a contracting rectangle $ R $, is established in Theorem \ref{prop-boundedness}. Specifically, the solution satisfies
$$
(v(x,t), u(x,t)) \in R, \quad \forall x \in \mathbb{R}, \quad t \geq 0.
$$

To derive the required estimate, we fix $ T > 0 $ and choose $ M $ sufficiently small such that if $ \|V\|_Y < M $, then
$$
2(1+a)\norm{\tilde{v}}_Y + 3\norm{\tilde{v}}_Y^2 < a - \frac{1}{\gamma}.
$$
Additionally, for a given $ \mu > 0 $, there exist constants $ C_1(\mu), C_2(\mu) $ satisfying:
\begin{enumerate}
    \item $ \|v\|_Y \leq C_1(\mu) $,
    \item $ \|\varphi_1\|_Y \leq C_2(\mu) $,
    \item $ \|\varphi_2\|_Y \leq 3C_1(\mu)+ a+1$.
\end{enumerate}

And, there exist constants $\varepsilon_1$ such that
\begin{enumerate}
    \item $\sup_{x \in \R} |\delta(x) - \sigma| < \varepsilon_1$,
\end{enumerate}
Under these conditions, the initial problem \eqref{linear-error-equation} admits a unique solution $ (F_v, F_w) \in C([0, T] ; \mathcal{B} \times \mathcal{B}) $, as guaranteed by Lemma \ref{lemma-linear-error}. Furthermore, Lemma \ref{lemma-nonlinear-error} ensures the existence of a unique solution $ (E_v, E_w) \in C([0, T] ; \mathcal{B} \times \mathcal{B}) $ to \eqref{difference-system}, which satisfies
$$
|E_v(x, t)| \leq \mu, \quad |E_w(x, t)| \leq \mu, \quad \forall x \in \mathbb{R}, \quad t \in [0, T].
$$

Consequently, we obtain the estimates
$$
| v - \tilde{v}| \leq \mu, \quad |u - \tilde{u}| \leq \mu, \quad \forall x \in \mathbb{R}, \quad t \in [0, T].
$$
Furthermore, as $T>0$ was arbitrary, the result holds for all $t>0$.

\bibliographystyle{acm}
\bibliography{bibliography}

@article{fhn,
title = {Impulses and Physiological States in Theoretical Models of Nerve Membrane},
journal = {Biophysical Journal},
volume = {1},
number = {6},
pages = {445-466},
year = {1961},
issn = {0006-3495},
doi = {https://doi.org/10.1016/S0006-3495(61)86902-6},
url = {https://www.sciencedirect.com/science/article/pii/S0006349561869026},
author = {Richard FitzHugh}
}

@ARTICLE{4066548nagumo,
  author={Nagumo, J. and Arimoto, S. and Yoshizawa, S.},
  journal={Proceedings of the IRE}, 
  title={An Active Pulse Transmission Line Simulating Nerve Axon}, 
  year={1962},
  volume={50},
  number={10},
  pages={2061-2070},
  doi={10.1109/JRPROC.1962.288235}
}

@article{rinzel1973traveling,
  title={Traveling wave solutions of a nerve conduction equation},
  author={Rinzel, John and Keller, Joseph B},
  journal={Biophysical journal},
  volume={13},
  number={12},
  pages={1313--1337},
  year={1973},
  publisher={Elsevier}
}

@article {MR889639KEENER,
    AUTHOR = {Keener, James P.},
     TITLE = {Propagation and its failure in coupled systems of discrete
              excitable cells},
   JOURNAL = {SIAM J. Appl. Math.},
  FJOURNAL = {SIAM Journal on Applied Mathematics},
    VOLUME = {47},
      YEAR = {1987},
    NUMBER = {3},
     PAGES = {556--572},
      ISSN = {0036-1399},
   MRCLASS = {92A09 (34G20 35R15)},
  MRNUMBER = {889639},
MRREVIEWER = {Harlan\ W.\ Stech},
       DOI = {10.1137/0147038},
       URL = {https://doi.org/10.1137/0147038},
}

@article {MR393759HASTING,
    AUTHOR = {Hastings, S. P.},
     TITLE = {On the existence of homoclinic and periodic orbits for the
              {F}itzhugh-{N}agumo equations},
   JOURNAL = {Quart. J. Math. Oxford Ser. (2)},
  FJOURNAL = {The Quarterly Journal of Mathematics. Oxford. Second Series},
    VOLUME = {27},
      YEAR = {1976},
    NUMBER = {105},
     PAGES = {123--134},
      ISSN = {0033-5606,1464-3847},
   MRCLASS = {35B10 (34C25)},
  MRNUMBER = {393759},
MRREVIEWER = {Joel\ Smoller},
       DOI = {10.1093/qmath/27.1.123},
       URL = {https://doi.org/10.1093/qmath/27.1.123},
}

@article{rall1969time,
  title={Time constants and electrotonic length of membrane cylinders and neurons},
  author={Rall, Wilfrid},
  journal={Biophysical journal},
  volume={9},
  number={12},
  pages={1483--1508},
  year={1969},
  publisher={Elsevier}
}

@article{kleber2004basic,
  title={Basic mechanisms of cardiac impulse propagation and associated arrhythmias},
  author={Kl{\'e}ber, Andr{\'e} G and Rudy, Yoram},
  journal={Physiological reviews},
  volume={84},
  number={2},
  pages={431--488},
  year={2004},
  publisher={American Physiological Society}
}

@article{clayton2011models,
  title={Models of cardiac tissue electrophysiology: progress, challenges and open questions},
  author={Clayton, RH and Bernus, Olivier and Cherry, EM and Dierckx, Hans and Fenton, Flavio H and Mirabella, Lucia and Panfilov, Alexander V and Sachse, Frank B and Seemann, Gunnar and Zhang, H},
  journal={Progress in biophysics and molecular biology},
  volume={104},
  number={1-3},
  pages={22--48},
  year={2011},
  publisher={Elsevier}
}

@article {MR3773759cornwell,
    AUTHOR = {Cornwell, Paul and Jones, Christopher K. R. T.},
     TITLE = {On the existence and stability of fast traveling waves in a
              doubly diffusive {F}itz{H}ugh-{N}agumo system},
   JOURNAL = {SIAM J. Appl. Dyn. Syst.},
  FJOURNAL = {SIAM Journal on Applied Dynamical Systems},
    VOLUME = {17},
      YEAR = {2018},
    NUMBER = {1},
     PAGES = {754--787},
      ISSN = {1536-0040},
   MRCLASS = {35K40 (35B25 35B35 35C07 35K57 35P15 37J05)},
  MRNUMBER = {3773759},
MRREVIEWER = {Paulo\ R.\ Zingano},
       DOI = {10.1137/17M1149432},
       URL = {https://doi.org/10.1137/17M1149432},
}

@article{cebrian2024sixgalens,
  title={Six decades of the FitzHugh--Nagumo model: A guide through its spatio-temporal dynamics and influence across disciplines},
  author={Cebri{\'a}n-Lacasa, Daniel and Parra-Rivas, Pedro and Ruiz-Reyn{\'e}s, Daniel and Gelens, Lendert},
  journal={Physics Reports},
  volume={1096},
  pages={1--39},
  year={2024},
  publisher={Elsevier}
}

@article {cerpaode,
    AUTHOR = {Cerpa, Eduardo and Courdurier, Mat\'ias and Hern\'andez,
              Esteban and Medina, Leonel E. and Paduro, Esteban},
     TITLE = {A partially averaged system to model neuron responses to
              interferential current stimulation},
   JOURNAL = {J. Math. Biol.},
  FJOURNAL = {Journal of Mathematical Biology},
    VOLUME = {86},
      YEAR = {2023},
    NUMBER = {1},
     PAGES = {Paper No. 8, 38},
      ISSN = {0303-6812,1432-1416},
   MRCLASS = {92-08 (37N25 92C30)},
  MRNUMBER = {4519609},
MRREVIEWER = {Jason\ M.\ Graham},
       DOI = {10.1007/s00285-022-01839-8},
       URL = {https://doi.org/10.1007/s00285-022-01839-8},
}

@misc{jhoncerpa,
      title={Analysis and stability of the FitzHugh-Nagumo system with rapidly oscillating Neumann boundary condition}, 
      author={Eduardo Cerpa and Alexander Rodriguez},
      year={under review},
      eprint={2305.00123},
      archivePrefix={arXiv},
      primaryClass={math.AP}
}

@article {cerpa2023approximation,
    AUTHOR = {Cerpa, Eduardo and Courdurier, Mat\'ias and Hern\'andez,
              Esteban and Medina, Leonel E. and Paduro, Esteban},
     TITLE = {Approximation and stability results for the parabolic
              {F}itz{H}ugh-{N}agumo system with combined rapidly oscillating
              sources},
   JOURNAL = {Discrete Contin. Dyn. Syst.},
  FJOURNAL = {Discrete and Continuous Dynamical Systems},
    VOLUME = {50},
      YEAR = {2026},
     PAGES = {69--102},
      ISSN = {1078-0947,1553-5231},
   MRCLASS = {35Q92 (35K55)},
  MRNUMBER = {5066925},
       DOI = {10.3934/dcds.2026001},
       URL = {https://doi-org.pucdechile.idm.oclc.org/10.3934/dcds.2026001},
}

@article{RAUCH197812,
title = {Qualitative theory of the FitzHugh-Nagumo equations},
journal = {Advances in Mathematics},
volume = {27},
number = {1},
pages = {12-44},
year = {1978},
issn = {0001-8708},
doi = {https://doi.org/10.1016/0001-8708(78)90075-0},
url = {https://www.sciencedirect.com/science/article/pii/0001870878900750},
author = {Jeffrey Rauch and Joel Smoller}
}

@article {SCHONBEK1978119,
    AUTHOR = {Schonbek, Maria Elena},
     TITLE = {Boundary value problems for the {F}itzhugh-{N}agumo equations},
   JOURNAL = {J. Differential Equations},
  FJOURNAL = {Journal of Differential Equations},
    VOLUME = {30},
      YEAR = {1978},
    NUMBER = {1},
     PAGES = {119--147},
      ISSN = {0022-0396,1090-2732},
   MRCLASS = {35K55 (35Q20 92A05)},
  MRNUMBER = {511479},
MRREVIEWER = {Kuang-Ho\ Chen},
       DOI = {10.1016/0022-0396(78)90027-X},
       URL = {https://doi.org/10.1016/0022-0396(78)90027-X},
}

@book {MR2263523,
    AUTHOR = {Izhikevich, Eugene M.},
     TITLE = {Dynamical systems in neuroscience: the geometry of
              excitability and bursting},
    SERIES = {Computational Neuroscience},
 PUBLISHER = {MIT Press, Cambridge, MA},
      YEAR = {2007},
     PAGES = {xvi+441},
      ISBN = {978-0-262-09043-8; 0-262-09043-0},
   MRCLASS = {37N25 (34C60 37C10 37G99 92C20)},
  MRNUMBER = {2263523},
MRREVIEWER = {Jonathan\ E.\ Rubin},
}

@article{hodgkin1952quantitative,
author = {Hodgkin, A. L. and Huxley, A. F.},
title = {A quantitative description of membrane current and its application to conduction and excitation in nerve},
journal = {The Journal of Physiology},
volume = {117},
number = {4},
pages = {500-544},
doi = {https://doi.org/10.1113/jphysiol.1952.sp004764},
url = {https://physoc.onlinelibrary.wiley.com/doi/abs/10.1113/jphysiol.1952.sp004764},
eprint = {https://physoc.onlinelibrary.wiley.com/doi/pdf/10.1113/jphysiol.1952.sp004764},
year = {1952}
}

@article{rattay1986analysis,
  title={Analysis of models for external stimulation of axons},
  author={Rattay, Frank},
  journal={IEEE transactions on biomedical engineering},
  number={10},
  pages={974--977},
  year={1986},
  publisher={IEEE}
}

@book {MR1952568,
    AUTHOR = {Murray, J. D.},
     TITLE = {Mathematical biology. {II}},
    SERIES = {Interdisciplinary Applied Mathematics},
    VOLUME = {18},
   EDITION = {Third},
      NOTE = {Spatial models and biomedical applications},
 PUBLISHER = {Springer-Verlag, New York},
      YEAR = {2003},
     PAGES = {xxvi+811},
      ISBN = {0-387-95228-4},
   MRCLASS = {92-02 (92B05 92C05 92D30)},
  MRNUMBER = {1952568},
MRREVIEWER = {Trachette\ L.\ Jackson},
}

@article {MR381744,
    AUTHOR = {Hastings, S. P.},
     TITLE = {Some mathematical problems from neurobiology},
   JOURNAL = {Amer. Math. Monthly},
  FJOURNAL = {American Mathematical Monthly},
    VOLUME = {82},
      YEAR = {1975},
    NUMBER = {9},
     PAGES = {881--895},
      ISSN = {0002-9890,1930-0972},
   MRCLASS = {92A05 (34C25)},
  MRNUMBER = {381744},
MRREVIEWER = {Brian\ D.\ Sleeman},
       DOI = {10.2307/2318490},
       URL = {https://doi.org/10.2307/2318490},
}

@article {MR445091,
    AUTHOR = {Conway, E. D. and Smoller, J. A.},
     TITLE = {A comparison technique for systems of reaction-diffusion
              equations},
   JOURNAL = {Comm. Partial Differential Equations},
  FJOURNAL = {Communications in Partial Differential Equations},
    VOLUME = {2},
      YEAR = {1977},
    NUMBER = {7},
     PAGES = {679--697},
      ISSN = {0360-5302,1532-4133},
   MRCLASS = {35B05 (92A05)},
  MRNUMBER = {445091},
MRREVIEWER = {Richard\ J.\ Field},
       DOI = {10.1080/03605307708820045},
       URL = {https://doi.org/10.1080/03605307708820045},
}

@article {MR1036209,
    AUTHOR = {Val\`encia, Marta},
     TITLE = {On invariant regions and asymptotic bounds for semilinear
              partial differential equations},
   JOURNAL = {Nonlinear Anal.},
  FJOURNAL = {Nonlinear Analysis. Theory, Methods \& Applications. An
              International Multidisciplinary Journal},
    VOLUME = {14},
      YEAR = {1990},
    NUMBER = {3},
     PAGES = {217--230},
      ISSN = {0362-546X,1873-5215},
   MRCLASS = {35K55 (35B35 35B45 35Q80 58F10)},
  MRNUMBER = {1036209},
MRREVIEWER = {H.\ J.\ Kuiper},
       DOI = {10.1016/0362-546X(90)90030-K},
       URL = {https://doi.org/10.1016/0362-546X(90)90030-K},
}

@article {MR3262177,
    AUTHOR = {Boukerrioua, K.},
     TITLE = {Global existence of solution for reaction diffusion systems
              with a full matrix of diffusion coefficients},
   JOURNAL = {Surv. Math. Appl.},
  FJOURNAL = {Surveys in Mathematics and its Applications},
    VOLUME = {9},
      YEAR = {2014},
     PAGES = {105--115},
      ISSN = {1843-7265,1842-6298},
   MRCLASS = {35K51 (35A01 35K57)},
  MRNUMBER = {3262177},
}

@article {MR425397,
    AUTHOR = {Rauch, Jeffrey},
     TITLE = {Global existence for the {F}itzhugh-{N}agumo equations},
   JOURNAL = {Comm. Partial Differential Equations},
  FJOURNAL = {Communications in Partial Differential Equations},
    VOLUME = {1},
      YEAR = {1976},
    NUMBER = {6},
     PAGES = {609--621},
      ISSN = {0360-5302,1532-4133},
   MRCLASS = {35Q99 (92A05)},
  MRNUMBER = {425397},
MRREVIEWER = {J.\ W.\ Jerome},
       DOI = {10.1080/03605307608820023},
       URL = {https://doi.org/10.1080/03605307608820023},
}

@article {MR3395129cartersandstede,
    AUTHOR = {Carter, Paul and Sandstede, Bj\"orn},
     TITLE = {Fast pulses with oscillatory tails in the
              {F}itz{H}ugh-{N}agumo system},
   JOURNAL = {SIAM J. Math. Anal.},
  FJOURNAL = {SIAM Journal on Mathematical Analysis},
    VOLUME = {47},
      YEAR = {2015},
    NUMBER = {5},
     PAGES = {3393--3441},
      ISSN = {0036-1410,1095-7154},
   MRCLASS = {35K58 (34C37 34E15 35B25 35B44 35C07 35M30)},
  MRNUMBER = {3395129},
MRREVIEWER = {Camelia\ Alexandra\ Pop},
       DOI = {10.1137/140999177},
       URL = {https://doi.org/10.1137/140999177},
}

@article {MR3551276carterrijk,
    AUTHOR = {Carter, Paul and de Rijk, Bj\"orn and Sandstede, Bj\"orn},
     TITLE = {Stability of traveling pulses with oscillatory tails in the
              {F}itz{H}ugh-{N}agumo system},
   JOURNAL = {J. Nonlinear Sci.},
  FJOURNAL = {Journal of Nonlinear Science},
    VOLUME = {26},
      YEAR = {2016},
    NUMBER = {5},
     PAGES = {1369--1444},
      ISSN = {0938-8974,1432-1467},
   MRCLASS = {35K57 (35B25 35B35 35C07 35M30 35P15)},
  MRNUMBER = {3551276},
MRREVIEWER = {Daniel\ \v Sev\v covi\v c},
       DOI = {10.1007/s00332-016-9308-7},
       URL = {https://doi.org/10.1007/s00332-016-9308-7},
}

@article {MR3927583,
    AUTHOR = {Mao, Yiqiu},
     TITLE = {Dynamic transitions of the {F}itzhugh-{N}agumo equations on a
              finite domain},
   JOURNAL = {Discrete Contin. Dyn. Syst. Ser. B},
  FJOURNAL = {Discrete and Continuous Dynamical Systems. Series B. A Journal
              Bridging Mathematics and Sciences},
    VOLUME = {23},
      YEAR = {2018},
    NUMBER = {9},
     PAGES = {3935--3947},
      ISSN = {1531-3492,1553-524X},
   MRCLASS = {92C20 (35K57 35M33 35Q92)},
  MRNUMBER = {3927583},
       DOI = {10.3934/dcdsb.2018118},
       URL = {https://doi.org/10.3934/dcdsb.2018118},
}

@article {MR4583896,
    AUTHOR = {Al-Juaifri, Ghassan A. and Harfash, Akil J.},
     TITLE = {Analysis of a nonlinear reaction-diffusion system of the
              {F}itz{H}ugh-{N}agumo type with {R}obin boundary conditions},
   JOURNAL = {Ric. Mat.},
  FJOURNAL = {Ricerche di Matematica. A Journal of Pure and Applied
              Mathematics},
    VOLUME = {72},
      YEAR = {2023},
    NUMBER = {1},
     PAGES = {335--357},
      ISSN = {0035-5038,1827-3491},
   MRCLASS = {35K57 (65M12)},
  MRNUMBER = {4583896},
MRREVIEWER = {Jaywan\ Chung},
       DOI = {10.1007/s11587-022-00711-7},
       URL = {https://doi.org/10.1007/s11587-022-00711-7},
}

@article {MR4887572rubio,
    AUTHOR = {Rubio-Mercedes, C. E. and Ver\~ao, Glauce Barbosa},
     TITLE = {Numerical simulation of traveling waves in the
              {F}itz{H}ugh-{N}agumo system via equilibria of nonlocal
              equations},
   JOURNAL = {Comput. Appl. Math.},
  FJOURNAL = {Computational \& Applied Mathematics},
    VOLUME = {44},
      YEAR = {2025},
    NUMBER = {5},
     PAGES = {Paper No. 218, 19},
      ISSN = {2238-3603,1807-0302},
   MRCLASS = {65J20 (35C07 35K57 65N35)},
  MRNUMBER = {4887572},
       DOI = {10.1007/s40314-025-03182-2},
       URL = {https://doi.org/10.1007/s40314-025-03182-2},
}

@article {MR4860996ammatouch,
    AUTHOR = {Amattouch, M. R.},
     TITLE = {A study of the {F}itz{H}ugh-{N}agumo model with diffusion
              term},
   JOURNAL = {Bol. Soc. Parana. Mat. (3)},
  FJOURNAL = {Boletim da Sociedade Paranaense de Matem\'atica. 3rd S\'erie},
    VOLUME = {43},
      YEAR = {2025},
     PAGES = {19},
      ISSN = {0037-8712,2175-1188},
   MRCLASS = {35B25 (26A33 35K57 65M06)},
  MRNUMBER = {4860996},
       DOI = {10.5269/bspm.70405},
       URL = {https://doi.org/10.5269/bspm.70405},
}

@article {MR4419132sanchez,
    AUTHOR = {S\'anchez-Gardu\~no, Faustino and Castellanos, V\'ictor},
     TITLE = {Traveling wave solutions for nonlinear reaction-diffusion
              equations as dynamical systems problems},
   JOURNAL = {Lobachevskii J. Math.},
  FJOURNAL = {Lobachevskii Journal of Mathematics},
    VOLUME = {43},
      YEAR = {2022},
    NUMBER = {1},
     PAGES = {141--161},
      ISSN = {1995-0802,1818-9962},
   MRCLASS = {37C29 (35C07 37D10)},
  MRNUMBER = {4419132},
MRREVIEWER = {Ronald\ Adams},
       DOI = {10.1134/s1995080222040199},
       URL = {https://doi.org/10.1134/s1995080222040199},
}

@article {MR3385151chen2,
    AUTHOR = {Chen, Chao-Nien and Choi, Y. S.},
     TITLE = {Traveling pulse solutions to {F}itz{H}ugh-{N}agumo equations},
   JOURNAL = {Calc. Var. Partial Differential Equations},
  FJOURNAL = {Calculus of Variations and Partial Differential Equations},
    VOLUME = {54},
      YEAR = {2015},
    NUMBER = {1},
     PAGES = {1--45},
      ISSN = {0944-2669,1432-0835},
   MRCLASS = {35J47 (34C37 35A24 35C07 35K40 35K57)},
  MRNUMBER = {3385151},
MRREVIEWER = {Peng\ Feng},
       DOI = {10.1007/s00526-014-0776-z},
       URL = {https://doi.org/10.1007/s00526-014-0776-z},
}

@article {MR3281845arioli,
    AUTHOR = {Arioli, Gianni and Koch, Hans},
     TITLE = {Existence and stability of traveling pulse solutions of the
              {F}itz{H}ugh-{N}agumo equation},
   JOURNAL = {Nonlinear Anal.},
  FJOURNAL = {Nonlinear Analysis. Theory, Methods \& Applications. An
              International Multidisciplinary Journal},
    VOLUME = {113},
      YEAR = {2015},
     PAGES = {51--70},
      ISSN = {0362-546X,1873-5215},
   MRCLASS = {92C20 (34D05 35B10 35B35 35C07 35K57 37L15)},
  MRNUMBER = {3281845},
MRREVIEWER = {Joseph\ L.\ Shomberg},
       DOI = {10.1016/j.na.2014.09.023},
       URL = {https://doi.org/10.1016/j.na.2014.09.023},
}

@article {MR4484075gambino,
    AUTHOR = {Gambino, Gaetana and Giunta, Valeria and Lombardo, Maria
              Carmela and Rubino, Gianfranco},
     TITLE = {Cross-diffusion effects on stationary pattern formation in the
              {F}itz{H}ugh-{N}agumo model},
   JOURNAL = {Discrete Contin. Dyn. Syst. Ser. B},
  FJOURNAL = {Discrete and Continuous Dynamical Systems. Series B. A Journal
              Bridging Mathematics and Sciences},
    VOLUME = {27},
      YEAR = {2022},
    NUMBER = {12},
     PAGES = {7783--7816},
      ISSN = {1531-3492,1553-524X},
   MRCLASS = {35K57 (35B36 35Q92)},
  MRNUMBER = {4484075},
       DOI = {10.3934/dcdsb.2022063},
       URL = {https://doi.org/10.3934/dcdsb.2022063},
}

@article {MR2494028faina,
    AUTHOR = {Berezovskaya, Faina and Camacho, Erika and Wirkus, Stephen and
              Karev, Georgy},
     TITLE = {``{T}raveling wave'' solutions of {F}itzhugh model with
              cross-diffusion},
   JOURNAL = {Math. Biosci. Eng.},
  FJOURNAL = {Mathematical Biosciences and Engineering. MBE},
    VOLUME = {5},
      YEAR = {2008},
    NUMBER = {2},
     PAGES = {239--260},
      ISSN = {1547-1063,1551-0018},
   MRCLASS = {92C20 (34C23 35Q80 37N25)},
  MRNUMBER = {2494028},
MRREVIEWER = {Carlo\ R.\ Laing},
       DOI = {10.3934/mbe.2008.5.239},
       URL = {https://doi.org/10.3934/mbe.2008.5.239},
}

@article {MR716853collins,
    AUTHOR = {Collins, Clyde},
     TITLE = {Diffusion dependence of the {F}itz{H}ugh-{N}agumo equations},
   JOURNAL = {Trans. Amer. Math. Soc.},
  FJOURNAL = {Transactions of the American Mathematical Society},
    VOLUME = {280},
      YEAR = {1983},
    NUMBER = {2},
     PAGES = {833--839},
      ISSN = {0002-9947,1088-6850},
   MRCLASS = {35K57 (35B99 92A09)},
  MRNUMBER = {716853},
MRREVIEWER = {Piero\ de Mottoni},
       DOI = {10.2307/1999649},
       URL = {https://doi.org/10.2307/1999649},
}

@article {MR4274265orkun,
    AUTHOR = {Tasbozan, Orkun},
     TITLE = {A popular reaction-diffusion model fractional
              {F}itzhugh-{N}agumo equation: analytical and numerical
              treatment},
   JOURNAL = {Appl. Math. J. Chinese Univ. Ser. B},
  FJOURNAL = {Applied Mathematics. A Journal of Chinese Universities. Ser.
              B},
    VOLUME = {36},
      YEAR = {2021},
    NUMBER = {2},
     PAGES = {218--228},
      ISSN = {1005-1031,1993-0445},
   MRCLASS = {35M10 (26A33 35C05 35R11 92C20)},
  MRNUMBER = {4274265},
       DOI = {10.1007/s11766-021-3810-x},
       URL = {https://doi.org/10.1007/s11766-021-3810-x},
}

@article {MR4308218cordoni,
    AUTHOR = {Cordoni, Francesco and Di Persio, Luca},
     TITLE = {Optimal control of the {F}itz{H}ugh-{N}agumo stochastic model
              with nonlinear diffusion},
   JOURNAL = {Appl. Math. Optim.},
  FJOURNAL = {Applied Mathematics and Optimization},
    VOLUME = {84},
      YEAR = {2021},
    NUMBER = {3},
     PAGES = {2947--2968},
      ISSN = {0095-4616,1432-0606},
   MRCLASS = {93E20 (47H06 49J55 60H15 91G80)},
  MRNUMBER = {4308218},
MRREVIEWER = {M.\ D.\ Marcozzi},
       DOI = {10.1007/s00245-020-09735-7},
       URL = {https://doi.org/10.1007/s00245-020-09735-7},
}

@book{friedman2008partial,
  title={Partial differential equations of parabolic type},
  author={Friedman, Avner},
  year={2008},
  publisher={Courier Dover Publications}
}

@article{pascucci2018parametrix,
  title={The parametrix method for parabolic SPDEs},
  author={Pascucci, Andrea and Pesce, Antonello},
  journal={Stochastic Processes and their Applications},
  volume={130},
  number={10},
  pages={6226--6245},
  year={2020},
  publisher={Elsevier}
}
\end{document}